\documentclass[12pt,oneside,reqno]{amsart}   	
\usepackage{geometry}
\usepackage{float}
\allowdisplaybreaks
\usepackage{graphicx}				
\usepackage{amssymb}
\usepackage{amsmath}
\usepackage[all]{xy}
\usepackage{mathtools}
\usepackage{tikz-cd}
\usepackage{enumitem}
\usepackage{tikz-cd}
\usepackage{appendix}

\usepackage[%
bookmarks=true,
colorlinks,
linkcolor=blue,
urlcolor=blue,
citecolor=blue,
plainpages=false,
pdfpagelabels,
final,
breaklinks=true\between
]{hyperref}
\usepackage{hyperref}
\usepackage[numbers,sort&compress]{natbib}

\newtheorem{thm}{Theorem}[section]
\newtheorem{lem}[thm]{Lemma}

\newtheorem{rem}{Remark}[section]

\numberwithin{equation}{section}

\def\Pb{\ifmmode{\Bbb P}\else{$\Bbb P$}\fi}
\def\Z{\ifmmode{\Bbb Z}\else{$\Bbb Z$}\fi}
\def\N{\ifmmode{\Bbb N}\else{$\Bbb N$}\fi}
\def\C{\ifmmode{\Bbb C}\else{$\Bbb C$}\fi}
\def\R{\ifmmode{\Bbb R}\else{$\Bbb R$}\fi}
\def\S{\ifmmode{S^2}\else{$S^2$}\fi}

\def\diam{\operatorname{diam}}

\def\S{\cal S}

\newenvironment{pf}{\paragraph{Proof:}}{\hfill$\square$ \newline}

\begin{document}

\title[rotationally symmetric self shrinkers]{Compactness and finiteness theorems for rotationally symmetric self shrinkers}
\author{Alexander Mramor}

\address{Department of Mathematics, Johns Hopkins University, Baltimore, MD, 21231}
\email{amramor1@jhu.edu}

\begin{abstract} In this note we first show a compactness theorem for rotationally symmetric self shrinkers of entropy less than 2, concluding that there are entropy minimizing self shrinkers diffeomorphic to $S^1 \times S^{n-1}$ for each $n \geq 2$ in the class of rotationally symmetric self shrinkers. Assuming extra symmetry, namely that the profile curve is convex, we remove the entropy assumption. Supposing the profile curve is additionally reflection symmetric we show there are only finitely many such shrinkers up to rigid motion.
\end{abstract}
\maketitle
\section{Introduction}

Self shrinkers  $M^n \subset \R^{n+1}$, that is surfaces satisfying 
\begin{equation}\label{se}
H - \frac{\langle x, \nu \rangle}{2} = 0
\end{equation} 
 are models for singularities of the mean curvature flow but outside of some convexity conditions or strict entropy bounds (c.f. \cite{BW1, BW2, H2, CM1}) they are far from completely understood. The most well understood case seems to be closed genus 0 self shrinkers in $\R^3$, which by work of Brendle \cite{Bren} must be the round sphere of radius $\sqrt{2}$. Almost nothing is known about general self shrinkers of more complicated topology, one of the few partial results, due to the author and S. Wang \cite{MW}, is that closed self shrinkers in $\R^3$ must be ``unknotted.'' 
$\medskip$

Perhaps the next natural question following Brendle then is what can be said about self shrinking tori in $\R^3$ or, in higher dimensions, self shrinking ``donuts'' -- hypersurfaces diffeomoprhic to $S^1 \times S^{n-1}$. For example, is the Angenent torus unique amongst embedded self shrinking donuts? The purpose of this note is to provide some compactness and discreteness results as evidence in answering this question amongst the class of rotationally symmetric self shrinkers and specializations thereof. 
$\medskip$

Given $G, \Lambda > 0$, Colding and Minicozzi \cite{CM, CM1} showed that the space of self shrinking surfaces of genus $g < G$ and entropy less than or equal to $\Lambda$ in $\R^3$ was compact in the $C^\infty_{loc}$ topology. The entropy assumption was later weakened to just assuming properness by Ding and Xin \cite{DX}. More recently, Sun and Z. Wang \cite{SW} refined this result to show that the space of self shrinkers of bounded entropy and \textit{fixed} genus $g$ in $\R^3$ was compact. Namely for every genus $g$, if there exists a self shrinker of bounded entropy of that genus there is a self shrinker of least entropy of that genus. In higher dimensions, up to $n = 6$, Barbossa, Sharp and Wei \cite{BSW} show that Colding and Minicozzi's result was true assuming an additional index bound. This note is entirely concerned with the codimension 1 case, but we also mention that in higher codimension Chen and Ma showed in \cite{CMa} some compactness results for Lagrangian self shrinkers in $\C^2$ (appropriately modifying the definition of self shrinker above). 
$\medskip$

The first result of this note, which like all the others in this note heavily exploits rotational symmetry, is a compactness theorem for closed rotationally symmetric self shrinkers which holds in every dimension: 
\begin{thm}\label{thm1} Let $\Lambda < 2$ and $n \geq 2$. Then the set $\Sigma(n, \Lambda)$ of closed rotationally symmetric embedded self shrinkers of entropy less than or equal to $\Lambda$ is compact in the $C^\infty_{loc}$ topology.
$\medskip$

 In particular for every $n \geq 2$, there is a rotationally symmetric embedded self shrinking torus of least entropy.
\end{thm}
Here entropy is in the sense of Colding and Minicozzi \cite{CM1} and for self shrinkers is simply the Gaussian area; the naturalness of the entropy quantity is described in section 2 below. Note that by work of Kleene and M{\o}ller and later Song \cite{KM, AS} that such self shrinkers are either the standard round shrinking spheres or embedded donuts. This set is also nontrivial (hence, the conclusion) by the construction of Drugan and Ngyuen \cite{DN} and, in $\R^3$, Angenent \cite{Ang} by way of the computation of Berchenko-Kogan \cite{YBK}\footnote{Angenent's construction produces shrinking donuts in all dimensions, but \cite{YBK} only concerns the surface case - as posed in question 6.1 therein presumably their analysis can be extended}. 
$\medskip$

In the general, nonrotationally symmetric case the entropy condition needed to prove related curvature bounds for self shrinkers is significantly more restrictive because the links of minimal cones which might be encountered (see the proof of lemma 3.1 below) is not generally well understood, and it is not clear if diameter bounds for closed self shrinkers are to be expected. 
$\medskip$

Next we impose extra symmetry; by assuming the profile curve is convex, we show the following compactness theorem which holds without any entropy bound:
\begin{thm}\label{thmconvex} Let $n \geq 3$ and denote by $\Sigma_{con}(n)$ the set of closed, embedded rotationally symmetric self shrinkers of $\R^{n+1}$ with convex profile curve. Then this set is compact in the $C^\infty_{loc}$ topology. 
\end{thm}

It seems plausible that the profile curves must in fact always be convex - to the authors knowledge there are no examples of rotationally symmetric embedded self shrinkers with profile curve which is known to not be convex.  Obviously this would be the case if the Angenent donuts are unique. 
$\medskip$

Lastly, we discuss the nature of the entropies attainable by compact spaces of rotationally symmetric self shrinkers and, assuming even more symmetry, we prove a finiteness theorem for them using the same technique (note that Angenent's torus is a member of the set in item (3) by its construction): 
\begin{thm}\label{thmref} Where $\Sigma(n, \Lambda)$ is as in theorem \ref{thm1}, and $\Sigma_{con}$ is as in theorem \ref{thmconvex}: 

\begin{enumerate}
 \item The entropy functional $\lambda$ attains only finitely many values on $\Sigma(n, \Lambda)$ for a fixed $n \geq 2$ and $\Lambda < 2$. 
 \item Similarly, $\lambda$ obtains only finitely many values on the set $\Sigma_{con}(n)$ for a fixed $n \geq 3$.
 \item Additionally, denote by $\Sigma_{refsym}(n)$ the set of rotationally symmetric embedded self shrinkers in $\Sigma_{con}$ whose profile curves are also reflection symmetric across the line perpendicular to their axis of rotation. Then for each $n \geq 3$ there are only finitely many elements in $\Sigma_{refsym}(n)$ up to rigid motion.
 
  \end{enumerate}
\end{thm} 

In particular item (1) and (2) seems to suggest that all self shrinking donuts in those sets have convex and reflection symmetric profile curve and hence have only finitely many elements, but there may certainly simply be Jacobi fields for the metric as well. 
$\medskip$

$\textbf{Acknowledgements.}$ The author is indebted to S. Angenent for generously sharing his ideas concerning the Poincar\'e map centrally used in the proof of theorem \ref{thmref} below. He is also grateful to J. Bernstein for some informative conversations on self shrinkers of low entropy and T. Bourni and M. Langford for their encouragement and interest in this work. Finally, he thanks the anonymous referee for their careful suggestions which helped improve the quality of this article. 

\section{Background on the mean curvature flow and justification of entropy} 

The mean curvature flow, where in most generality one deforms a hypersurface $x: M^n \to N^{n+k}$ by 
\begin{equation}
\frac{dx}{dt} = \vec{H}
\end{equation}
 is a parabolic system and enjoys many of the same properties as the classical heat equation. It is not strictly elliptic but under some assumptions of bounded geometry it is solvable for short times. When $M$ is a hypersurface (so $\vec{H} = - H \nu$, where $\nu$ is the outward unit normal) its flow satisfies the so-called comparison principle (also known as the avoidance principle), which says that two initially disjoint surfaces stay disjoint under the flow. 
 $\medskip$
 
  In particular, by comparison with a sufficiently large enveloping sphere, every compact hypersurface in Euclidean space must develop a singularity by some time $T_{s} < \infty$. In this paper our flows will will always be hypersurfaces of $\R^{n+1}$. We typically denote the flow associated to a given hypersurface $M$ by $M_t$. 
 $\medskip$
 
 To study the singularities that can arise, Huisken in \cite{H1} introduced the following: 
 \begin{equation}
\Phi_{x_0, t_0}(x,t) = \frac{1}{(4 \pi (t_0 - t))^{n/2}} \cdot \text{exp}\left ( - \frac{|x-x_0|^2}{4 (t_0 - t)} \right ), \text{ }t < t_0.
\end{equation} 
Then Huisken's montonicity (theorem 3.1 in \cite{H1}) says the integral of $\Phi_{x_0, t_0}$ is nonincreasing under the flow; more precisely 
\begin{thm}({Huisken monotonicity}) If $M_t$ is a surface flowing by the mean curvature flow for $t < t_0$, then we have the formula
\begin{equation}\label{variation} 
\frac{d}{dt} \int_{M_t} \Phi_{x_0, t_0}(x, t) d\mu_t = - \int_{M_t} \Phi_{x_0, t_0}(x,t) \left | H + \frac{\langle x, \nu \rangle }{2(t_0 - t)}\right |^2 d\mu_t
\end{equation}
\end{thm}

The right hand side readily implies that singularities are modeled, in a precise sense \footnote{namely when one performs a \textit{tangent flow} i.e. parabolically rescaling at a fixed point - certain other blowups such as Hamiltons type II blowup will often result in different models, such as translating solitons}, by surfaces satisfying the shrinker equation \ref{se} as discused by Huisken in \cite{H1} and under weaker assumptions by Ilmanen in \cite{I1}. 
$\medskip$

More geometrically, these surfaces correspond to the $t= -1$ timeslice of the (ancient) flow $\{\Sigma_t\}_{t\in(-\infty,0)}$ of the form $\Sigma_t=\sqrt{-t}\Sigma_{-1}$ - this flow is by dilations and hence is ``self similar.'' 
$\medskip$

 Huisken's monotonicty suggests that there should be a related metic for which equation \ref{variation} is roughly the first variation of area, so that self shrinkers are minimal surfaces in this metric. And indeed, self shrinkers are minimal surfaces in the Gaussian metric $G = e^{-|x|^2/4} g_{\mathrm{flat}}$. The area of a surface in this metric being given by the Gaussian area: 
\begin{equation}
A(M) := \int_M e^{-\frac{|x|^2}{4}} d\mu
\end{equation}
This notion of area is less than ideal for studying the flow for a couple reasons. Firstly the area funcitonal doesn't ``see'' regions far away from the origin in a sense, because the weight above decays exponentially -- a notion of area that is basepoint independent is much more useful in singularity analysis since ultimately one takes blowups to study singularities.
$\medskip$

 Secondly the metric turns out to be $f$-Ricci positive in the sense of Bakry and Emery and so has many of the same properties as Ricci positive metrics - see for example \cite{WW}. This is unfortunate, because in analogy to the Ricci positive case there are no stable minimal surfaces -- morally there should be such a class for the right notion of area because singularites as explained above are in many natural cases inevitable.  Then one should in principle be able to ``perturb away'' singularities modeled on unstable self shrinkers, and also so that they are very rare, in doing so simplify the study of the flow. We explore this in more detail. Indeed, the Jacobi operator for minimal surfaces in the Gaussian metric is given by:
\begin{equation}
L = \Delta + |A|^2  -\frac{1}{2} \langle x, \nabla(\cdot) \rangle + \frac{1}{2}
\end{equation} 
Because of the constant term on the RHS $L$ has nontrivial index in virtually all situations:
\begin{lem} Self shrinkers of polynomial volume growth are unstable in the Gaussian metric. 
\end{lem}

In their fundamental paper \cite{CM} Colding and Minicozzi introduced a very powerful new quantity called the entropy to study the mean curvature flow which deals with both issues at once and hence is arguably the ``right'' notion of area to pair with the mean curvature flow\footnote{See \cite{YL} for a discussion of entropy for higher codimension submanifolds of $\R^N$}. To elaborate, consider a hypersurface $\Sigma^k \subset \R^{\ell}$; then given $x_0 \in \R^{\ell}$ and $r > 0$ define the functional $F_{x_0, r}$ by 
\begin{equation}
F_{x_0, r}(\Sigma) = \frac{1}{(4 \pi t_0)^{k/2}} \int_\Sigma e^\frac{-|x - x_0|^2}{4r} d\mu
\end{equation} 
Colding and Minicozzi then define the entropy $\lambda(\Sigma)$ of a submanifold to be the supremum over all $F_{x_0, r}$ functionals:
\begin{equation} 
\lambda(\Sigma) = \sup\limits_{x_0, r} F_{x_0, r}(\Sigma) 
\end{equation} 
Huisken monotonicity in fact implies that the entropy is monotone under the flow. Note that equivalently $\lambda(\Sigma)$ is the supremum of $F_{0,1}$ when we vary over rescalings (changing $r$) and translations (choice of $x_0$) -- in particular there is no ``prefered point'' in $\R^{n+1}$ in defining the entropy. For hypersurfaces with polynomial growth this supremum is attained and, for self shrinkers $\Sigma$, $\lambda(\Sigma) = F_{0,1}(\Sigma)$ -- their area in the Gaussian metric. 
$\medskip$

Hence as promised self shrinkers are critical points for the entropy and this time there are stable such critical points, because we now consider variations in the $x$ and $r$ coordinates, which we now explain. If $\Sigma_s$ is a normal variation of a self shrinker $\Sigma$ and $x_s, r_s$ are variations with $x_0=0, r = 1$, 
\begin{equation}
\partial_s \mid_{s = 0} \Sigma_s = f\nu, \partial_s \mid_{s=0} x_s = y, \text{ and }\partial_s \mid_{s = 0} r_s = h
\end{equation} 
The second variation formula one finds for the entropy of $\Sigma$ is:
\begin{equation}\label{sv}
 (4 \pi)^{-n/2} \int_\Sigma (-fLf + 2fhH - h^2H^2 f \langle y, \nu \rangle - \frac{\langle y, \nu \rangle^2}{2})e^{\frac{-|x|^2}{4}} d\mu
 \end{equation}
 where $L$ is as before the Jacobi operator for the Gaussian metric.
 $\medskip$
 
One can easily check that, where $v$ is a constant vector field on $\R^n$, both $\langle v, \nu \rangle$ and $H$ are eigenfunctions with eigenvalues $-1, -\frac{1}{2}$ respectively for $L$; $LH = H$ and $L\langle v, \nu \rangle = \frac{1}{2}\langle v, \nu \rangle$. If $-1$ is not the lowest eigenvalue for $L$, then one could find a a function $f$, perpendincular (in the appropriate weighted $L^2$ space) to both $H$ and $\langle v, \nu \rangle$, such that $Lf = cf$ for $c > 0$. By this orthogonality one sees from the second variation formula \ref{sv} that $f$ can be used to build an entropy decreasing variation. So $-1$ must be the lowest eigenvalue, and so $H$ must have a sign. The mean convex self shrinkers are under very mild assumptions known to be generalized round cylinders $S^k \times \R^{n-k}$, so we have morally justified the following: 
\begin{thm}(Theorem 0.12 in \cite{CM}) Suppose that $\Sigma$ is a smooth complete embedded self-shrinker without boundary and with polynomial volume growth.
\begin{enumerate} 
\item If $\Sigma$ is not equal to $S^k \times \R^{n-k}$, then there is a graph $\widetilde{\Sigma}$ over $\Sigma$ of a function with arbitrarily small $C^m$ norm (for any fixed $m$) so that $\lambda(\widetilde{\Sigma}) < \lambda(\Sigma)$
\item If $\Sigma$ is not $S^n$ and does not split off a line, then the function in (1) can be taken to have compact support. 
\end{enumerate}
\end{thm} 
 Conversely, when $H$ does have a sign and is compact, the self shrinker is entropy stable (see lemma 4.23 in \cite{CM} -- $F$-stability is in most cases equivalent to entropy stability).
 $\medskip$

We see then that Colding and Minicozzi's notion of entropy is a very natural, and perhaps the best, notion of area for studying surfaces in the context of the mean curvature flow.

\section{Proof of theorem \ref{thm1}} 

Throughout $n$ and $\Lambda$ are implicitly fixed. First we prove the following curvature bound: 
\begin{lem}\label{curvbound} There is a function $f: \R_+ \to \R_+$ so that in the ball $B(0,R)$, $|A| \leq f(R)$ for $M \in \Sigma$. 
\end{lem}
\begin{pf}
To see this, suppose there is $R > 0$ and $M^i \in \Sigma$ with points $p^i \in B(0,R)$ such that $|A_{M^i}|(p^i) \to \infty$. 
By the self shrinker equation, we see that $H \sim R$ and hence is bounded in the ball $B(0,2R)$.  Thus rescaling by $|A|$ and recentering at $p^i$ then, we get subsequential convergence to a nonflat (since $|A| = 1$) minimal surface $N$ - note that if $p_i \to \partial B(0,R)$, it might be the case that $N$ has unbounded curvature. Our claim is that $N$ is a catenoid centered at the origin however. 
$\medskip$

To show this we argue that $\lim\limits_{i \to \infty} p^i = 0$. Indeed suppose that there was some subsequence (further relabeling to $p^i$) which was bounded away from the origin by some $\delta > 0$. Then by the rotational symmetry, since the curvatures on the surfaces at $p^i$ associated to rotation are on the order $\frac{1}{|p^i|}  < \frac{1}{\delta} < \infty$, $N$ must split off a (flat) $n-1$ plane. Since $N$ is minimal then the profile curve must itself have vanishing curvature, contradicting that $N$ is nonflat. 
$\medskip$

Thus $p^i \to 0$ and as a consequence $|A| \leq 1$ on $N$. Since $p^i \to 0$ and each of the $M^i$ are rotationally symmetric, it then follows that $N$ is a catenoid. The entropy of the catenoid is 2 however by the monotnicity formula; to see it is at least 2 note its blowdown is a puntured plane with multiplicity 2, which has itself entropy equal to 2, and entropy is lower semicontinuous. This of course contradicts that the entropy of each of the $M^i$ is less than $\Lambda < 2$. 
\end{pf}
$\medskip$

Now we prove a diameter bound - this is where the closedness assumption is used: 
\begin{lem}\label{diambound} Denote by $\Sigma_{closed}(n, \Lambda)$ the set of closed elements of $\Sigma$. Then there is a $D(n, \Lambda) < \infty$ for which all self shrinkers in $\Sigma_{closed}$ have diameter bounded by $D$. 
\end{lem} 
\begin{pf}
Again suppose that this wasn't true, and consider a sequence $M^i$ with diam$(M^i) \to \infty$. Denote by $p^i \in M^i$ points which achieve the farthest distance to the origin. At these points we see that the position vector $x$ and the unit normal of $M^i$ are parallel so that $H = |p^i|/2$. Of course, by the rotational symmetry these points correspond to a $n-1$ dimensional sphere of radius $|p^i|$ so these points aren't unique. We now study the geometry of $M^i$ at these points. 
$\medskip$

After possibly tilting the axis of symmetry, we may assume without loss of generality, because the group of rotations $O(n+1)$ of $\R^{n+1}$ is compact, that each of the $p_i$ are all in a fixed direction, say in the direction of the vector $e_{n+1} = (0,0, \ldots,1)$ along the $x_{n+1}$-axis.
$\medskip$

 Since $|p^i| \to \infty$, $\frac{1}{\sqrt{|p^i|}} \to 0$ so that the position vector on the balls $B(p^i, \frac{1}{\sqrt{|p^i|}})$ converges to $|p|e_{n+1}$, or more precisely that $|x - |p^i|e_{n+1}| = O(\frac{1}{\sqrt{|p^i|}})$\footnote{We will abuse notation by implicitly assuming the position vector $x$, etc. belongs to the $M_i$ for the same $i$ which appears in the given inequality} \footnote{Here when we use ``big O'' notation the asymptotic is taken in terms of $i$}. Hence in the balls $B(p^i, \frac{1}{\sqrt{|p^i|}})$ the self shrinker equation gives:
 \begin{equation}\label{r1}
H- |x| \frac{\langle e_{n+1}, \nu \rangle}{2} = O(\frac{1}{\sqrt{|p^i|}})
 \end{equation}
  Now we recall that the $n-1$ prcincipal curvatures $\lambda_k$ of $M^i$ at $x$ in the rotational directions have value $\lambda_k < 1/|x|$, so that $|H(x) - \kappa(x)| = O(\frac{1}{|p^i|})$ for $x \in B(p^i, \frac{1}{\sqrt{|p^i|}})$, where $\kappa(x)$ is the geodesic curvature of the profile curve at $x$, so we may refine \ref{r1} to say
  \begin{equation}\label{r2}
   \kappa - |x| \frac{\langle e_{n+1}, \nu \rangle}{2} = O(\frac{1}{\sqrt{|p^i|}})
 \end{equation}
for $x \in B(p^i, \frac{1}{\sqrt{|p^i|}})$. 
 $\medskip$
 
 Now we recenter the points $p^i$ to the origin and rescale by $|p^i|$, so that the balls $B(p^i, \frac{1}{\sqrt{|p^i|}})$ are mapped to the balls $B(0, \sqrt{|p^i|})$. We will get a sequence of surface $\widetilde{M^i}$, where naturally $\widetilde{x}$, $\widetilde{\kappa}$, $\widetilde{\nu}$ denote their position vector, curvature of profile curve, and outward normal respectively, and under this coordinate change \ref{r2} transforms to: 
   \begin{equation}\label{r3}
 |p^i| \widetilde{\kappa} - \left | \frac{\widetilde{x}}{|p^i|} +p^i \right | \frac{\langle e_{n+1}, \widetilde\nu \rangle}{2} = O(\frac{1}{\sqrt{|p^i|}})
 \end{equation}
for $\widetilde{x} \in B(p^i, \sqrt{|p^i|})$. Dividing through by $|p^i|$ we then get: 
   \begin{equation}\label{r3}
 \widetilde{\kappa} - \left | \frac{\widetilde{x}}{|p^i|^2} +\frac{p^i}{|p^i|} \right | \frac{\langle e_{n+1}, \widetilde\nu \rangle}{2} = O(\frac{1}{|p^i|^{3/2}})
 \end{equation}

  We wish to pass to the limit but first we need to check that the curvature on each of the $\widetilde{M^i}$ will be bounded. 
 $\medskip$
 
 For $i$ sufficiently large since the rotational principal curvatures in $B(p^i, \frac{1}{\sqrt{|p^i|}})$ are on the order of $\frac{1}{|p^i|}$ as already mentioned: 
\begin{equation}\label{r3}
|\kappa(x)| -1 < |H(x)| < |\kappa(x)| + 1, x \in B(p^i, \frac{1}{\sqrt{|p^i|}})
\end{equation} 
for $i$ sufficiently large. By that same reasoning, $\kappa$ in these balls is equivalent $|A|$ up to a bounded additive constant for $i$ sufficiently large. From the self shrinker equation, for closed self shrinkers $H$ attains its global maximum at the farthest points $p$ from the origin, where it attains value $|p|/2$, so that rescaling at $p^i$ as indicated above gives that $|\widetilde{A}|$ bounded by at least 3. 
$\medskip$

Passing to the limit then, we get a mean curvature flow which splits off an $n-1$ plane to give us a translating curve shortening flow translating at speed 1/2 in the $e_{n+1}$ direction; that is a curve shortening flow satisfying
  \begin{equation}\label{r4}
   \widetilde{\kappa} - \langle e_{n+1}/2, \widetilde{\nu} \rangle = 0
 \end{equation}
The only two translators in the plane are the grim reaper and the plane; and hence by nonflatness the limit profile curve is a grim reaper. This again yields a contradiction, since the entropy of the grim reaper is 2. 
\end{pf}

With these lemmas in place, let $M^i$ be a sequence of self shrinkers in the set $\Sigma(n, \Lambda)$ for a fixed choice of $n$ and $\Lambda < 2$. By lemma \ref{diambound} there exists $D < \infty$ for which the diameters of all the $M^i$ are bounded by $D$. By the standard covering and Arzela-Ascoli diagonal argument, from the local curvature bounds \ref{curvbound} we get locally smooth convergence of $M^i$ to a self shrinker $M$ which must also have diameter bounded by $D$; what remains to show is that $M$ is itself embedded. 
$\medskip$

If this were not the case, then for each $M_i$ there would be points $p^i$ and $q^i$ in a bounded set $D$, with normals perpendicular to each other, sot that $||p^i - q^i|| \to 0$ as $i \to \infty$. By the theorem, the self shrinkers $\Sigma^i$ are locally graphical with bounded gradient in uniform sized balls (independent of $i$) about each $p^i$ and $q^i$. It is easy to see then from using $F$-functionals centered at these points that $\lim_{i \to \infty} \lambda(M^i) \geq 2$, giving a contradiction. 
$\medskip$

Now we calculate as in \cite{BLM} that the rotationally symmetric shrinking donuts $T_n$ of Drugan and Nguyen have entropy less than 2 to justify the final conclusion of theorem \ref{thm1}: 
\begin{lem}\label{lem:entropy} The entropy of $T_n$ is strictly less than 2 for each $n$.  
\end{lem}
\begin{pf}
Denote by $L_n(C)$ the length of a curve $C$ in the plane induced by the Gaussian metric as explained in section 5 below, for which self shrinking donuts correspond to closed geodesics in the upper halfplane. The geodesics found in \cite{DN}, denoted $\gamma_\infty$, satisfy the following estimate:   
\begin{equation}
L_n(\gamma_\infty) < 2 \int\limits_0^\infty s^{n-1} e^{-s^2/4} ds = 2^n \Gamma\left(\frac{n}{2}\right).
\end{equation} 

The Gaussian area $A(T)$ of the corresponding self-shrinking $S^1\times S^{n-1}$ is then $\mathrm{Vol}(S^{n-1}) L_n(\gamma) < 2^n n \pi^{\frac{n}{2}}  \frac{\Gamma(\frac{n}{2})}{\Gamma(\frac{n}{2} + 1)}$. Now we recall from section 2 above that the entropy of a compact self-shrinker is equal to the $F$ functional $F_{0,1}$, which is simply the Gaussian area normalized so that the plane has value 1. Thus,
\begin{equation} 
\lambda(T) = F_{0,1}(T) = \frac{1}{(4\pi)^{n/2}}A(T) < \frac{2^n n \pi^{\frac{n}{2}}}{(4\pi)^{n/2}} \frac{\Gamma(\frac{n}{2})}{\Gamma(\frac{n}{2} + 1)} = n \frac{\Gamma(\frac{n}{2})}{\Gamma(\frac{n}{2} + 1)}  = 2
\end{equation} 
Which is the bound we claimed. 
\end{pf}
$\medskip$

Lemma \ref{lem:entropy} implies that for each $n$, there is a $\Lambda_n < 2$ for which the set $\Sigma(n ,\Lambda)$ contains self shrinkers diffeomorphic to $S^{n-1} \times S^1$ provided $\Lambda > \Lambda_n$ (alternatively as mentioned in the introduction, when $n =2$ the computation of Berchenko-Kogan \cite{YBK} gives the Angenent torus in $\R^3$ has entropy approximately $1.85122 <2$). Fixing $n$ and such a $\Lambda < 2$, the set $\Sigma(n, \Lambda)$ is compact; denote by $\lambda(n)$ the infimum of entropy over self shrinking donuts in $\Sigma$. Letting $D_k$ be a sequence of self shrinking donuts in $\Sigma$ so that $\lambda(D_k) < \lambda(n) + 1/k$, we may take a converging subsequence by compactness and since the convergence is in the smooth topology the limit must itself be a donut, giving the last assertion of theorem \ref{thm1}.

\section{Proof of Theorem \ref{thmconvex}}

Throughout $n \geq 3$ is fixed. The proof of theorem \ref{thmconvex} again comes down to proving curvature and diameter bounds but the techniques are a bit different and in some respects closer to the techniques in the previous works mentioned in the introduction. There are two cases to rule out:
\begin{enumerate}
\item There is a sequence $M^i \subset \Sigma_{con}$ with $ \lim\limits_{i \to \infty} \sup\limits_{p \in M^i} |A|(p) \to \infty$ but diameter uniformly bounded by some $D < \infty$. 
\item There is a sequence $M^i \subset \Sigma_{con}$ with $diam(M^i) \to \infty$. 
\end{enumerate}
Indeed, if case (2) is ruled out then there must be a diameter bound; ruling out case (1) then tells us that there must then be curvature bounds on elements of $\Sigma_{con}$ which will give its compactness essentially as before. 
$\medskip$

First we deal with case (1). As in the lemma \ref{curvbound} above, by recentering the points $p^i$ where $\sup\limits_{p \in M^i} |A|(p)$ is obtained to the origin and rescaling by $|A|(p^i)$ we get in the limit a catenoid $C$. Since $n \geq 3$ and $|A| = 1$ somewhere on the catenoid $C$ is contained in a slab of width independent of $M^i$, so that by the convexity assumption each of the $M^i$ are contained in slabs $S^i$ where their width $w_i \to 0$ as $i \to \infty$. 
$\medskip$

Now we use that the diameters of the $M^i$ are all uniformly bounded by some constant $D < \infty$. For $\overline{w}$ small enough, one can arrange that a grim reaper $G_t$ with width less than $\overline{w}$ will travel distance greater than $2D$ in time less than $1/2$. For $i$ large enough of course $w_i < \overline{w}$ so that, by arranging a grim reaper plane $G_t \times \R^{n-1}$ of width less than $\overline{w}$ engulfing $M^i$ with that the distance from the origin to the tip of the grim reaper is less than $2D$, we must have by the comparison principle that the $M^i_t$ (with the flow starting at $t = -1$) for $i$ large enough must develop a singularity before $t = -1/2$. This is a contradiction though since this occurs at $t = 0$ for self shrinkers. 
$\medskip$

We borrow some of the notation from case (1) to discuss case (2). Also without loss of generality, again since $O(n+1)$ is compact, we may suppose that the axis of rotation of each of the $M^i$ is $x_{n+1}$. First we discuss the limiting behavior if case (2) were to occur: 
$\medskip$

\begin{lem} As $i \to \infty$, $M^i \to W := \{\vec{x} \mid x_{n+1} = 0\} \setminus \{0\}$ with multiplicity 2 in the $C^\infty_{loc}$ topology. 
\end{lem} 
\begin{pf} 

It is clear by the convexity of the profile curves that a subsequential limit of the $M^i$ exists weakly and will be a cone with link given by a round (by rotational symmetry) $S^{n-1}$ sphere $S \subset S^n$; it will suffice to show that $S$ is a great equator. Because $diam(M^i) \to \infty$, since the profile curves are convex, and since the tips of the self shrinkers at extremal points are modeled on grim reapers as discussed in the proof of theorem \ref{thm1}, the width $w_i$ of the profile curves $P^i$ tend to 0 as $i \to \infty$. 
$\medskip$

Of course because $\diam(M^i) \to \infty$, we must have that for $i$ sufficiently large that all of the $M^i$ are donuts, since otherwise they would have to be the round shrinking sphere. Hence the $M^i$ after rescaling converge to the catenoid perpendicular to $e_1$ about the origin, implying by convexity of the profile curve that the link of the limit is an equator, giving the statement. 

\end{pf}

 Now we make an application of Simon's method \cite{Si} following Colding and Minicozzi \cite{CM} (perhaps the prototype compactness theorem using this method is due to Choi and Schoen \cite{CS}) to show that this implies the plane is stable in the Gaussian metric and get a contradiction; we sketch the details for the reader. 
$\medskip$

We continue to suppose the axis of symmetry of each of the $M^i$ is $x_{n+1}$. Away from the origin and by the convexity of the profile curves, we may write the $M^i$ as graphs of functions $f^i_+$ and $f^i_{-}$ over the plane $P =\{\vec{x} \mid x_{n+1} = 0 \}$. Since each of the $M^i$ are shrinkers the difference $f^i = f^i_+ - f^i_{-}$ satisfies $Lf^i = 0$ up to some correction terms which we can scale away: fixing some $y$ away from the origin and setting $g^i = \frac{f^i}{f^i(y)}$ by elliptic estimates we may let $i \to \infty$ to get a positive function $g$ over $P\setminus \{0\}$ such that $Lg = 0$ and $g(y) = 1$. 
$\medskip$

By using a local minimal foliation \cite{W1} of controlled geometry due to White along with the maximum principle one can then show that the origin is a removable singularity of $g$ and so $g$ extends to a positive (hence nonzero) function on $P$ such that $Lg = 0$. By Barta's lemma (see \cite{FCS} for Barta's lemma in the noncompact case) it then follows that the first eigenvalue of $L$ on any compact region of the plane is strictly positive, contradicting the instablity of the plane as a minimal surface in the Gaussian metric. 
$\medskip$

With cases (1) and (2) above ruled out we then get as before uniform curvature and diameter bounds for elements of $\Sigma_{con}$. What remains to show is that the limit of a sequence of elements $M^i \subset \Sigma_{con}$ also belongs there. Again we may suppose $M^i$ are all donuts; as before the point is to rule out multisheeted convergence in the limit. Instead of using an entropy condition to rule this out as in the proof of Theorem \ref{thm1}, this is now a consequence of the convexity of each of the profile curves -- if there is multisheeted convergence the curvature along $M^i$ must blowup, contradicting case (1).

\section{Proof of Theorem \ref{thmref}}

 We begin with defining some notations to discuss the analytic properties of profile curves of self shrinkers more carefully, following Angenent's notation in \cite{Ang} closely although there are some small but important differences which will be pointed out. Also, in the proof below we show items (2) and (3) first and then return to (1) at the end. 
 $\medskip$
 
 One can parameterize a rotationally symmetric surface locally as follows. Consider an immersed curve $\alpha(s): (a,b) \to \R^2$ given by $s \to (x(s), r(s))$ in the $xr$ plane, and by associating the $xr$ plane with the $x_1  x_2$ plane in $\R^{n+1}$ the rotation of the curve $\alpha$ about the $x$ axis is parameterized by $X(\omega, s) : S^{n-1} \times (a,b) = x(s)e_1 + r(s) \omega$, where $\omega$ are polar coordinates on the $n-1$ sphere into $(x_2, \ldots, x_{n+1}) = \R^n \subset \R^{n+1}$. The Gaussian metric in this case can then be reduced to a related metric on the $xr$ plane with line element given by 
\begin{equation}
(ds)^2 = r^{2(n-1)}e^{-(x^2 + r^2)/4} \{ (dx)^2 + (dr)^2 \}
\end{equation} 
Closed rotationally symmetric shrinking donuts then correspond to closed geodesics of the upper half plane in this metric. Given any point $(x,r)$ in the upper halfplane and angle $\theta \in (-\pi, 0)$, there exists a unique geodesic $(x,r)$ whose tangent at that point is $(\cos{(\theta)}, \sin{(\theta)})$. Along a geodesic we then have the following system of ODE: 
\begin{equation}\label{geoev}
x' = \cos{(\theta)}, r' =  \sin{(\theta)}, \theta' = \frac{x}{2} \sin{(\theta)} + ( \frac{n-1}{r} - \frac{r}{2}) \cos{(\theta)}
\end{equation}
Note that by the Cauchy-Kovalevskaya theorem that solutions to the system of ODE above must be analytic in $t$ since RHS are all analytic. Generalizing slightly Angenent's notation, let $(x_{R,\theta}(t), r_{R,\theta}(t), \theta_{R,\theta}(t)) = \Gamma_{R, \theta}: [0, T(R)) \to \R^3$ be the maximal solution of system \ref{geoev} with initial value $(0,R,\theta)$, let $T^* = T^*(R, \theta) > 0$ be the second time at which $x_R = 0$ occurs.
$\medskip$

 Note that time could concievably never occur for some initial data but by analytic dependence of analytic ODE on initial conditions (extend to the complex domain and use Theorem 1.1 in \cite{IY}) to see $\Gamma_{R, \theta}$ is analytic (and hence continuous)) in $R$ and $\theta$ so that it is at least well defined in neighborhoods of points where such an intersection does occur; this will suffice for our uses. When $T*$ is well defined, consider the Poincar\'e map $P: \R^2 \to \R^2$ defined by
\begin{equation}\label{pm}
P(R, \theta) = (r_{R, \theta}(T^*), \theta_{R, \theta}(T^*))
\end{equation}
Clearly closed geodesics which bound a convex region must be fixed points of $P$ (of course, there might be nonconvex closed geodesics which are fixed points of $P$ as well). By analytic dependence of analytic ODE on initial conditions $T^*$ is an analytic function in the initial data. Since $r_{R, \theta}$ and $\theta_{R, \theta}$ are analytic, $P$ is an analytic map on its domain of defintion (depending on when $T^*$ is defined). An important feature of it is the following: 
\begin{lem}\label{fpset} The fixed points of the Poincar\'e map $P: \R^2 \to \R^2$ defined above are either isolated or analytic curves in the plane.
\end{lem}
\begin{pf}
This follows from the analytic implicit function theorem applied to the function $P - I$, where $I$ is the identity map -- fixed points of $P$ then correspond to zeroes of $P-I$. If at a zero $\vec{z}$ $d(P-I)$ has rank 2, then the zero is isolated. Otherwise $dP$ has an eigenvector of eigenvalue 1, which by change of basis we may suppose without loss of generality is the vector $e_1$. Since $P$ is not identically equal to the identity map $dP - I$ must have rank precisely equal to 1 then, so we may then use the analytic implict function theorem to find an analytic function $g(x)$ so the zero set locally about $\vec{z}$ is given by the analytic curve $x \to (x, g(x))$. 
\end{pf} 

With this observation we now show part (2) of theorem \ref{thmref}. To start, consider a one parameter family of fixed points to the Poincar\'e map and denote by $N^s$ the corresponding family of self shrinkers and that by the mean value theorem, $\lambda(M^s) = F_{0,1}(N^s)$ for all $s$ is independent of $s$. 
$\medskip$

With this in mind, suppose that there were infinitely many values of the entropy obtained in any of the compact sets $\Sigma$ or $\Sigma_{con}$ (again, for an allowable choice of $\Lambda$ and $n$), and let $M^i$ be a sequence of self shrinkers in these sets such that $\lambda(M^i) \neq \lambda(M^j)$ for $i \neq j$. By the compactness, we may extract a converging subsequence and by lemma \ref{fpset} we see for $i$ large enough that each of the $M^i$ lay in a one parameter family. This is a contradiction from the observation in the previous paragraph however. 
$\medskip$

Now we move on to part (3). Note that the set $\Sigma_{refsym}$ is compact by theorem \ref{thmconvex} and the additional fact that the limit of reflection symmetric profile curves will itself be reflection symmetric. Now we consider the ``forgetful'' Poincar\'e map $P_f: \R \to \R$ given by $P_f(R) = \pi_1(P(R, 0)) = r_R(T^*(R,0))$. While some of its fixed points might not correspond to self shrinkers since the angle is ``forgotten,'' true reflection symmetric closed geodesics must give rise to fixed points of this map. 
$\medskip$

This again is an analytic map, and since it is not identically the identity map it must have isolated zeroes. This with the compactness gives (3) of theorem \ref{thmref}. 
$\medskip$

\begin{rem} Note that the reflection symmetry above was used because the initial angle of the geodesic ray is predetermined so that the ray is perpendicular to the line of symmetry; for any fixed choice of initial angle a similar statement is true.
\end{rem}

Now we finally discuss (1); it follows from the above by some simple modifications. Since it requires no extra machinery for the sake of exposition we first show only countably many values may be obtained by the entropy and then show the full statement in (1). Comparing with (2), note the convexity assumption can be dropped and replaced with the entropy condiiton in Theorem \ref{thm1} along with the stipulation the profile curve intersects the $r$ axis twice. Generally speaking by uniqueness of ODE a profile curve of a closed self shrinker can intersect the $r$-axis only transversely, so that by entropy arguments as in the proof of theorem \ref{thm1} the profile curve of a given self shrinker in $\Sigma(n, \Lambda)$ can intersect the $r$-axis only finitely many times. Now consider the filtration of $\Sigma(n, \Lambda)$ by the sets (each individually compact):
\begin{equation}
\Sigma_m(n, \Lambda) := \{ M \in \Sigma(n, \Lambda) \mid \text{the profile curve of }M\text{ intersects the }r\text{-axis }m\text{ times}\}
\end{equation} 
Following the arugment for item (2), $\lambda$ attains only finitely many values on $\Sigma_m$ for a fixed $m$ by, instead of using $T^*$ above in the defintion of the Poincar\'e map, one uses $T_m^*$ which stands for the $m$ time the profile curve intersects the $x$ axis. Since from the discussion above $\Sigma = \cup_{m \geq 2} \Sigma_m$, $\lambda$ must only obtain countably many values on $\Sigma(n, \Lambda)$. 
$\medskip$

To see there are actually only finitely many values of entropy obtained on the set $\Sigma(n, \Lambda)$ then, from the above it suffices to show that all but finitely many of the sets $\Sigma_m(n, \Lambda)$ are actually empty. In fact, it turns out that only $\Sigma_2(n, \Lambda)$ will be nonempty, by item (1) of Theorem 4 in \cite{KM} concluding the proof of theorem \ref{thmref}.

\end{document}